\documentclass{amsart}

\def\g{\gamma}

\def\p{\varphi}

\def\ra{\rightarrow}

\def\D{\partial}
\def\bx{{\bf x}}

\def\bxi{{\bf \xi}}
\def\l{\lambda}

\def\ep{\epsilon}
\def\a{\alpha}

\def\d{\delta}
\def\be{\begin{equation}}
\def\ee{\end{equation}}
\def\bea{\begin{eqnarray}}
\def\eea{\end{eqnarray}}
\def\ba{\begin{array}}
\def\ea{\end{array}}

\newtheorem{theorem}{Theorem}[section]
\newtheorem{lemma}[theorem]{Lemma}

\theoremstyle{definition}

\theoremstyle{remark}
\newtheorem{remark}[theorem]{Remark}

\numberwithin{equation}{section}

\newfont{\Bb}{msbm8 scaled\magstep{1}}
\newcommand{\rc}{\mbox{\Bb R}}

\usepackage{graphics}

\begin{document}
In the book "Differential Equations and Applications", 
vol.4,
Nova Sci. Publishers, New York, 2004

\title
[MRC for periodic structures]
{Modified Rayleigh Conjecture for scattering by periodic structures
}

\author{Alexander G. RAMM}
\address{ Department of Mathematics\\ Kansas State University\\
Manhattan, Kansas 66506-2602, USA}
\email{ramm@math.ksu.edu}

\author{Semion GUTMAN}
\address{Department of Mathematics\\ University of Oklahoma\\ Norman,
 OK 73019, USA}
\email{sgutman@ou.edu}


\begin{abstract} 
This paper contains a self-contained brief presentation of the scattering theory for 
periodic structures. Its main result is a theorem (the Modified Rayleigh 
Conjecture, or MRC), which gives a rigorous foundation for a numerical method
for solving the direct scattering problem for periodic structures. A
numerical example illustrating the procedure is presented.
\end{abstract}

\maketitle



\section{Introduction}
For simplicity we consider a 2-D setting, but our arguments can be as
easily  applied to $n$-dimensional problems, $n\geq 2$. 
Let $f : \rc\ra\rc,\ f(x+L)=f(x)$
be an $L$-periodic Lipschitz continuous function, and let $D$ be the
domain
\[
D=\{(x,y)\ :\ y\geq f(x),\ x\in \rc\}.
\]

Without loss of generality we assume that $f\geq 0$. If it is not, one can 
choose the origin so that this assumption is satisfied, because
$M:=\sup_{0\leq x \leq L} |f(x)|<\infty$.

Let $\bx=(x,y)$ and $u(\bx)$ be the total field satisfying
\begin{equation}\label{s1_1}
(\Delta +k^2)u=0,\quad \bx\in D, \quad k=const>0
\end{equation}
\begin{equation}\label{s1_2}
u=0 \quad \text{on} \quad S:\,=\D D,
\end{equation}
\begin{equation}\label{s1_3}
u=u_0+v,\quad u_0:\,=e^{ik\a\cdot\bx},
\end{equation}
where the unit vector $\a=(\cos\theta,-\sin\theta),\ 0<\theta<\pi/2$,
and $v(\bx)$ is the scattered field, whose asymptotic behavior as
$y\ra\infty$ will be specified below, and 
\begin{equation}\label{s1_4}
u(x+L,y)=\nu u(x,y),\quad u_x(x+L,y)=\nu u_x(x,y)\;\text{in}\; D,\quad
\nu:\,=e^{ikL\cos\theta}\,.
\end{equation}

Conditions (\ref{s1_4}) are the $qp$ ({\bf quasiperiodicity}) conditions. 
To find the proper radiation condition for the scattered field $v(\bx)$
consider the spectral problem
\begin{equation}\label{s1_5}
\p''+\l^2\p=0, \quad 0<x<L,
\end{equation}
\begin{equation}\label{s1_6}
\p(L)=\nu\p(0), \quad \p'(L)=\nu\p'(0)
\end{equation}
arising from the separation of variables in (\ref{s1_1})-(\ref{s1_4}).
This problem has a discrete spectrum, and its eigenfunctions form a
basis in $L^2(0,L)$. One has
\[
\p=Ae^{i\l x}+Be^{-i\l x},\quad A,B=const,
\]
\[
Ae^{i\l L}+Be^{-i\l L}=\nu(A+B),\quad i\l Ae^{i\l L}-i\l Be^{-i\l L}=i\l
\nu(A-B).
\]
Thus
\[
\begin{vmatrix}
e^{i\l L}-\nu & e^{-i\l L}-\nu \\
i\l(e^{i\l L}-\nu) & -i\l(e^{-i\l L}-\nu)
\end{vmatrix}
=0.
\]
So, $\ i\l(e^{i\l L}-\nu)(e^{-i\l L}-\nu)=0$. 
If $\l=0$, then $\p=A+Bx,\ A+BL=\nu A,\ B=\nu B$.  Since
$\nu=e^{ikL\cos\theta}$, one has no eigenvalue $\l =0$ unless
$kL\sin\theta=2\pi m,\ m>0$ is an integer. Let us assume that 
$kL\cos\theta\not=2\pi m$. Then
\[
e^{i\l L}=e^{ikL\cos\theta}\quad \text{or} \quad e^{-i\l
L}=e^{ikL\cos\theta}\,,
\] 
that is
\[
\l_j^+=k\cos\theta+\frac{2\pi j}L,\quad\text{or}\quad 
\l_j^-=-k\cos\theta+\frac{2\pi j}L,\quad j=0,\pm1,\pm2,\dots
\]
The corresponding eigenfunctions are $e^{i\l_j^+x}$ and $e^{-i\l_j^-x}$.
We will use the system $e^{i\l_j^+x}$, which forms an orthogonal bais 
in $L^2(0,L)$. One has:
\[
\int_0^L e^{i\l_j^+x}e^{-i\l_m^+x}\ dx=\int_0^Le^{\frac{2\pi i}{L}(j-m)}\
dx=0,\quad j\not=m.
\]
The normalized eigenfunctions are
\[
\p_j(x)=\frac{e^{i\l_j^+x}}{\sqrt{L}},\quad j=0,\pm1,\pm2,\dots
\]
These functions form an orthonormal basis of $L^2(0,L)$.
Let us look for $v(\bx)=v(x,y)$ of the form
\begin{equation}\label{s1_7}
v(x,y)=\sum_{j=-\infty}^{\infty} c_jv_j(y)\p_j(x), \quad y>M,\quad c_j=const.
\end{equation}
For $y>M$, equation (\ref{s1_1}) implies
\begin{equation}\label{s1_8}
v_j''+(k^2-\l_j^2)v_j=0.
\end{equation}
Let us assume that $\l_j^2\not=k^2$ for all $j$.
Then
\begin{equation}\label{s1_9}
v_j(y)=e^{i\mu_jy},
\end{equation}
where, for finitely many $j$, the set of which is denoted by $J$, one has:
\begin{equation}\label{s1_10}
\mu_j=(k^2-\l_j^2)^{1/2}>0,\quad\text{if}\quad \l_j^2<k^2, \,\, j\in J,
\end{equation}
and 
\begin{equation}\label{s1_11}
\mu_j=i(\l_j^2-k^2)^{1/2},\quad\text{if}\quad \l_j^2>k^2, \,\, j\notin J.
\end{equation}

The {\bf radiation condition} at infinity requires that the scattered
field $v(x,y)$ be representable in the form (\ref{s1_7}) with $v_j(y)$
defined by (\ref{s1_9})-(\ref{s1_11}).

The {\bf Periodic Scattering Problem} consists of finding the solution to 
(\ref{s1_1})-(\ref{s1_4}) satisfying the radiation condition (\ref{s1_7}), 
(\ref{s1_9})-(\ref{s1_11}).

The existence and uniqueness for such a scattering problem is
established in Section 2. Our presentation is essentially 
self-contained. In \cite{alber} the scattering by a periodic structure was 
considered earlier, and was based on a uniqueness theorem from \cite{e}.
 Our proofs differ from the proofs in  \cite{alber}. There are many papers
on scattering by periodic structures, of which we mention a few
\cite{alber}, \cite{albertsen}, \cite{bon},
\cite{bonram}, \cite{christiansen}, \cite{kazan1},\cite{kazan2},
\cite{millar1}, \cite{millar2}, \cite{petit},
\cite{ral1}. The Rayleigh conjecture is discussed in several of the above 
papers. It was shown (e.g.  \cite{petit},  \cite{baran})  that this 
conjecture is incorrect, in general. 
The modified Rayleigh conjecture is a theorem proved in \cite{r430} for
scattering by bounded obstacles.
A numerical method for solving obstacle scattering problems, based on the
modified Rayleigh conjecture is developed in \cite{gutmanramm1}.
The main results of our paper are: the modified Rayleigh
conjecture for periodic structures (Theorem 4.4) and a rigorous numerical 
method
for solving scattering problems by periodic structures, based on the
modified Rayleigh conjecture (Section 4). The proof of the limiting 
absorption principle (LAP) and the rigorous and self-contained development 
of the plane wave scattering theory by periodic structures is also of 
interest for broad audience. This theory is based partly on the ideas 
developed in \cite{rammb1}, \cite{r3}, \cite{r409}, \cite{r370}.
The proof of the key lemma 2.2 is based on a version of
Ramm's identity (2.16).
Numerical implementation of the method
for solving scattering problems by periodic structures, based on the
modified Rayleigh conjecture, is constructed using the 
approach developed 
in \cite{gutmanramm1} and in \cite{r475}. Applications to 
inverse 
problems are discussed in [18] and \cite{r470}.

\section{Periodic Scattering Problem}
Existence and uniqueness of solutions of the Periodic Scattering Problem
can be proved easily, if one establishes first the existence and
uniqueness of the resolvent kernel $G(x,y,\xi,\eta,k)$ of the Dirichlet
Laplacian in $D$:
\begin{equation}\label{s2_1}
(\Delta+k^2)G(x,y,\xi,\eta,k)=-\d(x-\xi)\d(y-\eta),\quad G=0\quad\text{on}
\quad S,
\end{equation}
\begin{equation}\label{s2_2}
G(x+L,y,\xi,\eta,k)=\nu G(x,y,\xi,\eta,k),\quad
G(x,y,\xi +L,\eta,k)=\overline{\nu} G(x,y,\xi,\eta,k),
\end{equation}
\begin{equation}\label{s2_3}
G_x(x+L,y,\xi,\eta,k)=\nu G_x(x,y,\xi,\eta,k),\quad
G_x(x,y,\xi +L,\eta,k)=\overline{\nu} G_x(x,y,\xi,\eta,k),
\end{equation}
and $G$ satisfies the LAP, see (2.5) below. 
The overbar here and below stands for the complex conjugation.
 
Indeed, if such a function $G$ exists, then $v$ can be found by the
Green's formula
\begin{equation}\label{s2_4}
v(x,y)=-\int_{S_L}u_0(\xi,\eta)G_N(x,y,\xi,\eta,k)\ ds,
\end{equation}
where $N$ is the unit normal vector to $S$ pointing into $D$.

To prove the existence and uniqueness of $G(x,y,\xi,\eta,k)$ define 
\[
\ell_0=-\Delta
\]
to be the Laplacian on the set of $C^2(D)$ quasiperiodic functions
vanishing on the boundary $S$, and vanishing near infinity. Let
\[
D_L:\,=\{(x,y)\ :\ 0\leq x\leq L,\quad (x,y)\in D\}.
\]
Then $D_L$ is a section of $D$, and $\ell_0$ is a symmetric operator in $L^2(D_L)$.
This operator is nonnegative , and therefore \cite{k} there exists its 
unique
selfadjoint Friedrichs' extension, which will be denoted by $\ell$.

Let $Im(k^2)>0$. Then there exists a unique resolvent operator $(\ell-k^2)^{-1}$. Thus its kernel
 $G(x,y,\xi,\eta,k)$ also exists and it is unique. To establish the
existence and uniqueness of the kernel for $k>0$ we are going to
prove the following 

{\bf Limiting Absorption Principle (LAP)}. Let $k>0,\ \epsilon>0$
and assume that $k^2$ is not equal to $\lambda_j^2$. Then 
the limit
\begin{equation}\label{s2_5}
\lim_{\epsilon\ra 0+}G(x,y,\xi,\eta,k+i\epsilon)=G(x,y,\xi,\eta,k),
\end{equation}
exists for all $(x,y)\in D,\ x\not=y$.
The proof is based on the following two lemmas.

\begin{lemma}\label{s2_l1} Let $\ 0<\ep<1$, and $a>2$. Then
\begin{equation}\label{s2_6}
\int_{D_L}\frac{|G(x,y,\xi,\eta,k+i\ep)|^2}{(1+\xi^2+
\eta^2)^{a/2}}\ d\xi d\eta\leq c, \end{equation} where
$c=const>0$ does not depend on $\ep$, and $(x,y)$ is running
on compact sets. 
\end{lemma}

\begin{proof}[Proof of Lemma \ref{s2_l1}]
It is sufficient to prove that the solution to the problem
\begin{equation}\label{s2_7}
(\Delta +k^2+i\ep)w_\ep=F,\ \text{in}\  D_L,
\ w_\ep\in L^2(D_L),\  w_\ep=0
\ \text{on}\ S_L
\end{equation}
\begin{equation}\label{s2_8}
w_\ep(x+L,y)=\nu w_\ep(x,y),\quad w_{\ep x}(x+L,y)=\nu w_{\ep x}(x,y),
\end{equation}
satisfies the estimate
\begin{equation}\label{s2_9}
N_\ep^2:\,=\sup_{0<\ep<1}\int_{D_L}\frac{|w_\ep(x,y)|^2}{(1+x^2+y^2)^{a/2}}\ dx dy:\,=N^2(w_\ep)\leq c,
\end{equation}
where $F\in C^\infty_0(D_L)$ is arbitrary, and $c=const>0$ is independent
of $\ep>0$.

If (\ref{s2_9}) fails, then $N_{\ep_n}\ra\infty,\ \ep_n\ra 0$. Define
$\psi_\ep:\,=w_\ep/N_\ep$, where $\ep:\,=\ep_n$. Then
$N(\psi_\ep)=1,\ \psi_\ep$ solves (\ref{s2_7}) (with $F$ replaced by
$F_\ep:\,=F/N_\ep$), and satisfies (\ref{s2_8}). From $N(\psi_\ep)=1$ it
follows that $\psi_\ep\rightharpoonup\psi$ as $\ep\ra 0$, where
$\rightharpoonup$ denotes the weak convergence in $L^2(D_L,1/(1+x^2+y^2)^{a/2}):\,=L^2_a$.
By elliptic estimates, $\psi_\ep\rightharpoonup\psi$ in
$H^2_{loc}(D_L)$, and therefore $\psi_\ep\ra\psi$ in
$L^2_{loc}(D_L)$. This and (\ref{s2_7})-(\ref{s2_8}) imply 
$\psi_\ep\ra\psi$ in $H^2_{loc}(D_L)$. Thus $\psi$ solves the
homogeneous ($F=0$) problem (\ref{s2_7})-(\ref{s2_8}).
If we prove that $\psi=0$,
then we get a contradiction, which shows that (\ref{s2_9}) holds. The
contradiction comes from the relationship
$0=N(\psi)=\lim_{\ep\ra 0}N(\psi_\ep)=1$. One proves that 
\begin{equation}\label{s2_star}
\lim_{\ep\ra 0}N(\psi_\ep)=N(\psi)
\end{equation}
as follows. If 
\[
(x,y)\in D_R:\,=\{(x,y)\ :\ f(x)\leq y\leq R,\ 0\leq x\leq L \},
\]
where $R>M$ is an arbitrary large fixed number, then $\lim_{\ep\ra
0}N(\psi_\ep \eta_R)=N(\psi \eta_R)$,
where
\[
\eta_R:\,=\begin{cases}
1, & f(x)<y<R,\\
0, & y>R.
\end{cases}
\]
In the region $D'_R=\{(x,y)\ :\ y>R,\ 0\leq x\leq L \}$, one has
$|\psi_\ep(x,y)|\leq c,\ (x,y)\in D'_R$. Thus
\[
\sup_{0<\ep<1}N(\psi_\ep(\chi_L-\eta_R))\leq O\left(\frac 1{R^\gamma}\right),\
0<\gamma<a-2.
\]
The desired result (\ref{s2_star}) follows.

To complete the proof let us show that the problem (\ref{s2_7})-(\ref{s2_8}), 
with $F=0$, and $\ep=0$, has only the trivial solution $w$, 
provided that $w$ is "outgoing" in the sense
\[
w_{jy}-i\mu_jw_j=o(1),\ \text{as}\ y\ra\infty,\ w_j:\,=\int_0^L
w\overline{\p_j}\ dx.
\]
One has
\begin{equation}\label{s2_10}
\lim_{R\ra\infty}\int_{S_R}(w\overline{w_y}-w_y\overline{w})\ ds=0,
\end{equation}
where $S_R:\,=\{(x,y)\ :\ y=R,\ 0\leq x\leq L \},\ ds=dx$ is the element
of the arclength of $S_R$, and the overbar stands for the complex
conjugate.

Let us outline the steps of the further argument.

{\it Step 1}: we prove that (\ref{s2_10}) implies
\begin{equation}\label{s2_11}
w\in L^2(D_L),\ |w|+|\nabla w|\leq ce^{-\g|y|}\,,\ \g=const>0,
\end{equation}
if $w$ is outgoing. 

{\it Step 2}: we prove that if $w\in L^2(D_L)$ solves (2.7)-(2.8),
with $F=\ep=0$, then $w=0$. 
Then we conclude that (2.9) (and (2.6)) holds, and, therefore, (2.5) 
holds. 

Let us prove (2.12). One has
\begin{equation}\label{s2_12}
\begin{split}
0 & =\int_{D_{LR}} [\bar w(\Delta+k^2)w-w(\Delta+k^2)\bar w]\ dxdy    \\
& =-\int_{S_L}(\bar w w_N-w\bar w_N )\ ds+\int_{S_R}(\bar w w_N-w\bar
w_N )\ ds \\
& = \int_{S_R}(\bar w w_N-w\bar w_N )\ ds,
\end{split}
\end{equation}
where the Dirichlet condition (\ref{s2_7}) was used, and the integrals 
over
the lines $x=0$ and $x=L$ are cancelled due to the $qp$ conditions
(\ref{s2_8}):
\[
\begin{split}
& \int_{x=0} (-\bar w w_x+w\bar w_x )\ dy+\int_{x=L}(\bar w w_x-w\bar
w_x )\ dy \\
& =\int_{x=0} ( w\bar w_x-\bar w w_x )\ dy-\int_{x=0}\nu\bar \nu(w\bar  w_x-\bar w
w_x )\ dy=0.
\end{split}
\]
Here we have used the relation $\nu\bar \nu=1$. Thus (\ref{s2_12}) implies
\begin{equation}\label{s2_13}
0=\int_{S_R}(\bar w w_y-w\bar
w_y )\ dx,\ \forall R>M. 
\end{equation}
If $w$ is outgoing, then (\ref{s2_13}) implies $w_j(y)=0$ for $j\in J$,
and $|w_j(y)|\leq e^{-\g|y|}\,,\ \g=const>0$, so (\ref{s2_11}) holds.
\end{proof}

\begin{lemma}\label{s2_l2}
Assume that $w\in L^2(D_L)$,  $w$ solves (\ref{s2_7}) with $\ep=0$ and 
$F=0$, and $w$ satisfies (\ref{s2_8}). Then $w=0$.
\end{lemma}

\begin{proof}[Proof of Lemma \ref{s2_l2}]
If $w$ solves equation (\ref{s2_7}) with $\ep=0$ and $F=0$, then 
$w=\sum_{j}w_j(y)\p_j(x)$. Since $\{\p_j(x)\}$ is an orthonormal basis
and $w\in L^2(D_L)$, it follows that $w_j(y)=0$ for all $j\in J$, and (\ref{s2_11})
holds. Let us use a version of Ramm's identity (\cite{r370}, p. 92),
which is valid for
any solution $w$ of equation (\ref{s1_1}) which is outgoing in the sense
that
\begin{equation}\label{s2_14}
w=\sum_{j}c_jv_j(y)\p_j(x),\ c_j=const,\ j\not\in J.
\end{equation}
Note, that $v_j(y)=\overline{v_j(y)}$ for $j\not\in J$.
The identity is:
\begin{equation}\label{s2_15}
0=(x_2 \bar w_{,2}w_{,j})_{,j}+\frac{(k^2|w|^2 x_2-|\nabla w|^2 x_2)_{,2}}{2}
+\frac{|\nabla w|^2-k^2|w|^2}2-|w_{,2}|^2,
\end{equation}
where $w_{,j}:\,=\D w/\D x_j,\ j=1,2,\ x_1=x,\ x_2=y$, over the repeated
indices one sums up, $|w|^2:\,=w\bar w$. The right-hand side of
(\ref{s2_15}) equals to
\[
\frac 12 [x_2(\bar w_{,2j}w_{,j}- w_{,2j}\bar w_{,j})+k^2x_2
( w_{,2}\bar w- \bar w_{,2} w)]=0,
\]
because $w_{,2}\bar w= \bar w_{,2} w$ for outgoing $w$.

One has
\begin{equation}\label{s2_16}
 |w|+|\nabla w|\leq ce^{-\g|y|}\,,\ \g=const>0,\ c=const>0.
\end{equation}
Let $R>\max f(x)$. Integrate (\ref{s2_15}) over $D_{LR}:\,=
\{(x,y)\ :\ (x,y)\in D_L,\ y\leq R\}$ and use Green's formula to get:
\begin{equation}\label{s2_17}
\begin{split}
0 & =-\lim_{R\to \infty} \int_{S_L\cup S_R}[x_2 
\bar w_{,2}w_{,j}N_j+\frac{(k^2|w|^2 x_2-
|\nabla w|^2x_2)N_2}2]\ ds\\
& -\lim_{R\to \infty}\int_{D_{LR}}|w_{,2}|^2\ dx_1dx_2,
\end{split}
\end{equation}
where $N$ is the normal pointing into $D_{LR}$, and we have used the
relation
\begin{equation}\label{s2_18}
\lim_{R\to \infty}\int_{D_{LR}}|\nabla w|^2\ dx_1dx_2=k^2
\lim_{R\to \infty}\int_{D_{LR}}|w|^2\ dx_1dx_2,
\end{equation}
which follows from the equation $\Delta w+k^2 w=0$, boundary condition
$w=0$ on $S$, quasiperiodicity of $w$, and from (\ref{s2_16}).
We have also used the relation $\bar w_{,2}w_{,j}N_j=x_2|\nabla 
w|^2N_2$, which follows from the condition $u=0$ on $S$.
 From
(\ref{s2_17}) one gets:
\begin{equation}\label{s2_19}
\lim_{R\to \infty}\int_{D_{LR}}|w_{,2}|^2\ dx_1dx_2=
-\frac 12 \int_{S_L}x_2N_2|\nabla w|^2\ ds.
\end{equation}
Since $f(x)$ is a graph, one has $N_2x_2\geq 0,$ and it follows from 
(\ref{s2_19}) that
$w_{,2}=0$, so $w=const$, and $const=0$ because $w|_S=0$.
Lemma \ref{s2_l2} is proved.
\end{proof}
\begin{remark}
Condition of the type
\begin{equation}\label{s2_20}
N_2x_2\geq 0\ \text{on}\ S_L
\end{equation}
was also used  in \cite{r370}.

The proof of Lemma 2.2 {\it is not valid
if the Neumann boundary condition is imposed on $S$.} 
\end{remark}

\section{Integral equations method}
In this Section we present another proof of the existence and uniqueness of 
the resolvent kernel $G$.
We want to construct a scattering theory quite similar 
to the one for the exterior of a bounded obstacle \cite{rammb1}. The {\it 
first step} is to construct an 
analog to the 
half-space Dirichlet Green's function. The function $g=g(\bx,\xi,k)$ can 
be
constructed analytically ($\bx=(x_1,x_2), \bxi=(\xi_1,\xi_2)$):
\begin{equation}\label{s3_01}
g(\bx,\bxi)=\sum_j\p_j(x_1)\overline{\p_j(\xi_1)}g_j(x_2,\xi_2,k),
\end{equation}
\[
g_j:=g_j(x_2,\xi_2,k)=\begin{cases}
v_j(x_2)\psi_j(\xi_2), & \ x_2>\xi_2\\
v_j(\xi_2)\psi_j(x_2), & \ x_2<\xi_2
\end{cases}
\]
\[
\psi_j=(\mu_j)^{-1}e^{i\mu_jb}\sin[\mu_j(\xi_2+b)],\, \, 
\mu_j=[k^2-\lambda_j^2]^{1/2},\quad v_j(x_2)=e^{i\mu_j x_2},
\]
where
\[
\psi_j''+(k^2-\l_j^2)\psi_j=0, \ \psi_j(-b)=0,\ W[v_j,\psi_j]=1,\
\lambda_j=k\cos (\theta)+\frac {2\pi j}{L},
\]
and $W[v,\psi]$ is the Wronskian. 

The function $g$ is analytic with respect to $k$ on the complex plain with 
cuts 
along the rays $\lambda_j-i\tau, \, 0\leq \tau <\infty, j=0, \pm 1, \pm 
2,...$, in particular, in the region $\Im k>0,$ up to the real
positive half-axis except for the set $\{\lambda_j\}_{j=0, \pm 1, \pm 2, 
...}$.

{\it Choose $b>0$ such that $k^2>0$ is not an eigenvalue of the problem:}
\begin{equation}\label{s3_1}
(\Delta +k^2)\psi=0,\quad \text{in}\ D_{-b}:\,=\{(x,y)\ :\ -b\leq y\leq f(x), \quad 0\leq x\leq L\}.
\end{equation}
\begin{equation}\label{s3_2}
\begin{split}
 & \psi|_{y=-b}=0, \quad \psi_N=0\ \text{on}\ S,\\
& \psi(x+L,y)=\nu \psi(x,y),\quad \psi_x(x+L,y)=\nu \psi_x(x,y).
\end{split} 
\end{equation}
One has
\begin{equation}\label{s3_3}
\begin{split}
& (\Delta +k^2)g=-\d(\bx-\bxi),\ \bx=(x_1,x_2),\ \bxi=(\xi_1,\xi_2),\\
& \bx\in \{(x,y)\ :\ -b< y<\infty, \quad 0\leq x\leq L\},
\end{split}
\end{equation}
\begin{equation}\label{s3_4}
g|_{y=-b}=0, 
\end{equation}
and
\begin{equation}\label{s3_5}
(\Delta +k^2)G=-\d(\bx-\bxi),\\ 
G=0\ \text{on}\ S,
\end{equation}
$G$ satisfies the $qp$ condition and the radiation
condition ( it is outgoing at infinity).

Multiply (\ref{s3_3}) by $G$, (\ref{s3_5}) by $g$, subtract
from the second equation the first one, integrate
over $D_{LR}$, and take $R\ra\infty$, to get
\begin{equation}\label{s3_6}
G=g+\int_{S_L}(Gg_N-G_Ng)ds=g-\int_{S_L}g\mu\ ds,\ \mu:\,=G_N|_{S_L}.
\end{equation}
The $qp$ condition allows one to cancel the integrals over the 
lateral boundary ($x=0$ and $x=L$), and the radiation condition allows one 
to have 
$$\lim_{R\to \infty} \int_{S_R}(Gg_N-G_Ng)ds=0.$$

Differentiate (\ref{s3_6}) to get 
\begin{equation}\label{s3_7}
\mu=-A\mu+2\frac{\D g}{\D N}\ \text{on}\ S_L,\ A\mu:\,=
2\int_{S_L}\frac{\D g(s,\sigma)}{\D N_s}\mu(\sigma)\ d\sigma.
\end{equation}
This is a Fredholm equation for $\mu$ in $L^2(S_L)$, if $S_L$ is
$C^{1,m},\ m>0$. The homogeneous equation (\ref{s3_7}) has only the
trivial solution: if $\mu+A\mu=0$, then the function
$\psi:\,=\int_{S_L}g\mu\ ds$ satisfies $\psi^+_N|_{S_L}=0$,
where
$\psi_N^+(\psi_N^-)$ is the normal derivative of $\psi$ from $D_{-b}
(D_L)$, and we use the known formula for the normal derivative of the
single layer potential at the boundary. The $\psi$
satisfies also  (\ref{s3_1}) and (\ref{s3_2}), and, by the choice 
of $b$, one has $\psi=0$ in $D_{-b}$. Also $\psi=0$ in $D_L$, because
$(\Delta+k^2)\psi=0$ in $D_L,\ \psi|_{S_L}=0$ (by the continuity of the
single layer potential), $\psi$ satisfies the $qp$ condition (because
$g$ satisfies it), and $\psi$ is outgoing (because $g$ is).

Since $\psi=0$ in $D_{-b}$ and in $D_L$, one concludes that 
$\mu=\psi_N^+-\psi_N^-$, where 
$\psi_N^+(\psi_N^-)$ is the normal derivative of $\psi$ from $D_{-b}
(D_L)$, and we use the jump relation for the normal derivative of the
single layer potential. 

Thus, we have proved the existence and
uniqueness of $\mu$, and, therefore, of  $G$, and got a representation 
formula
\begin{equation}\label{s3_8}
G=g-\int_{S_L}g\mu\ ds.
\end{equation}
This representation shows that the rate of decay of $G$ as $y\ra\infty$
is essentially the same as that of $g$.

The $G$ is analytic with respect to $k$ on the complex plain with cuts
along the rays $\lambda_j-i\tau, \, 0\leq \tau <\infty, j=0, \pm 1, \pm
2,...$, in particular, in the region $\Im k>0,$ up to the real
positive half-axis except for the set $\{\lambda_j\}_{j=0, \pm 1, \pm 2,
...}$. This follows from (\ref{s3_7}), (\ref{s3_8}), and the general
result \cite{rammb1}, p. 57, \cite{r189}, concerning analyticity of the solution to a 
Fredholm equation with respect to a parameter.

Suppose a bounded obstacle $D_0$ is placed inside $D_L$, $u=0$ on
$S_0=\partial D_0$, $S_0$ is a Lipschitz boundary. If $qp$ condition is 
imposed, then Green's function $G_0$ in the presence of the obstacle 
satisfies equations similar to (\ref{s3_8}) and (\ref{s3_7}):
\begin{equation}\label{s3_9}
G_0(x,y)=G(x,y)-\int_{S_0}G(x,s)\mu_0(s,y)\ ds, \quad \mu_0=G_{0N},
\end{equation}
where $N$ is the unit normal to $S_0$ pointing into $D_L$, and
\begin{equation}\label{s3_10}
\mu_0=-A_0\mu_0+2\frac{\D G}{\D N}\ \text{on}\ S_0,\ A_0\mu_0:\,=
2\int_{S_0}\frac{\D G(s,\sigma)}{\D N_s}\mu_0(\sigma)\ d\sigma.
\end{equation}
This is a Fredholm equation (with index zero). If $k^2$
is not an eigenvalue of the Neumann Laplacian in $D_0$ (=not exceptional),
then equation (\ref{s3_10}) is uniqueley solvable and, by (\ref{s3_9}), $G_0$ exists and 
is unique for this $k>0$. It is not known what are nontrivial sufficient 
conditions for $k>0$ to be not exceptional. The exceptional $k$ form a 
discrete countable set on the positive semi-axis $k>0$.
If the Neumann boundary condition is imposed on $S_L$, then, even in the 
absence of the obstacle $D_0$, it is not known if LAP holds, because the 
proof of Lemma 2.2 {\it is not valid for the Neumann boundary condition on 
$S_L$}.

 \section{Modified Rayleigh Conjecture 
(MRC)}
Rayleigh conjectured \cite{ral1} ("Rayleigh hypothesis") that the series
(\ref{s1_7}) converges up to the boundary $S_L$. This conjecture is
wrong (\cite{petit}) for some $f(x)$. Since the Rayleigh
hypothesis has been widely used for numerical solution of the scattering
problem by physicists and engineers, and because these practitioners
reported high instability of the numerical solution, and there are no
error estimates, we propose a modification of the Rayleigh conjecture,
which is a Theorem. This MRC (Modified Rayleigh Conjecture) can be used
for a numerical solution of the scattering problem, and it gives an error
estimate for this solution. Our arguments are very similar to the ones
in \cite{r430}. 

Rewrite the scattering problem (\ref{s1_1})-(\ref{s1_4}) as 
\begin{equation}\label{s4_1}
(\Delta+k^2)v=0\ \text{in}\ D,\ v=-u_0\ \text{on}\ S_L,
\end{equation}
where $v$ satisfies (\ref{s1_4}), and $v$ has representation
(\ref{s1_7}), that is, $v$ is "outgoing", it satisfies the radiation
condition. Fix an arbitrarily small $\ep>0$, and assume that
\begin{equation}\label{s4_2}
\|u_0+\sum_{|j|\leq j(\ep)}c_j(\ep)v_j(y)\p_j(x)\|\leq\ep,
\ 0\leq x\leq L,\ y=f(x),
\end{equation}
where $\|\cdot\|=\|\cdot\|_{L^2(S_L)}$.

\begin{lemma}\label{s4_l1} For any $\ep>0$, however small,
and for any $u_0\in L^2(S_L)$, there exists $j(\ep)$ and
$c_j(\ep)$ such that (\ref{s4_2}) holds. \end{lemma}
\begin{proof} Lemma \ref{s4_l1} follows from the
completeness of the system 
$\{\p_j(x)v_j(f(x))\}_{j=0, \pm
1, \pm2,....} $ in $L^2(S_L)$. Let us prove this
completeness. Assume that there is an $h\in L^2(S_L),\
h\not\equiv 0$
 such that
\begin{equation}\label{s4_3}
\int_{S_L}h\overline{\p_j(x)}v_j(f(x))\ ds=0
\end{equation}
for any $j$.
From (\ref{s4_3}) one derives (cf. \cite{rammb1}, p.162-163)
\begin{equation}\label{s4_4}
\psi(\bx):\,=\int_{S_L}hg(\bx,\bxi) d\xi=0,\ \bx\in D_{-b}.
\end{equation}
Thus $\psi=0$ in $D_L$, and $h=\psi_N^+-\psi_N^-=0$. Lemma \ref{s4_l1}
is proved.
\end{proof}

\begin{lemma}\label{s4_l2}
If (\ref{s4_2}) holds, then
\[
\||v(\bx)-\sum_{|j|\leq j(\ep)}c_j(\ep)v_j(y)\p_j(x)\||\leq c\ep
,\ 
\forall x,y\in D_L,
\ 0\leq x\leq L,\, y\geq f(x),
\]
where $c=const>0$ does not depend on $\ep, x, y,$
and $R$; $R>M$ is an arbitrary 
fixed number, and 
$\||w\||=\sup_{\bx \in D\setminus 
D_{LR}}|w(\bx)|+||w||_{H^{1/2}(D_{LR})}$.
\end{lemma}
\begin{proof}
Let $w:\,=v-\sum_{|j|\leq j(\ep)}c_j(\ep)v_j(y)\p_j(x)$. Then $w$ solves 
equation (\ref{s1_1}),
$w$ satisfies (\ref{s1_4}), $w$ is outgoing, and $\|w\|_{L^2(S_L)}\leq\ep$. 
One has (cf. (\ref{s2_4}))
\begin{equation}\label{s4_5}
w(\bx)=-\int_{S_L}wG_N(\bx,\bxi)\ ds.
\end{equation}
Thus (\ref{s4_2}), i.e. $\|w\|:\,=\|w\|_{L^2(S_L)}\leq\ep$, implies
\begin{equation}\label{s4_6}
|w(\bx)|_{y=R}\leq\|w\|_{L^2(S_L)}\|G_N(\bx,\bxi)\|_{L^2(S_L)}\leq c\ep,\ 
c=const>0,
\end{equation}
where $c$ is independent of $\ep$, and $R>\max f(x)$ is arbitrary.
Now let us use the elliptic inequality
\begin{equation}\label{s4_7}
\|w\|_{H^m(D_{LR})}\leq c\left( \|w\|_{H^{m-0.5}(S_L)} +  
\|w\|_{H^{m-0.5}(S_R)}     \right),
\end{equation}
where we have used the equation $\Delta w+k^2w=0$, and 
assumed that $k^2$ is not a Dirichlet 
eigenvalue of the Laplacian in $D_{LR}$, which can be done without loss
of generality, because one can vary $R$. 
The integer $m\geq 0$ is
arbitrary if $S_L$ is sufficiently smooth, and $m\leq 1$ if
 $S_L$ is Lipschitz. Taking 
$m=0.5$ and using (\ref{s4_2})
and (\ref{s4_6}) one gets 
\begin{equation}\label{s4_8}
\|w\|_{H^{1/2}(D_{LR})}\leq c\ep.
\end{equation}
Thus, in a neighborhood of $S_L$, we have proved estimate (\ref{s4_8}),
and in a complement of this neighborhood in $D_L$ we have 
 proved estimate (\ref{s4_6}). 
 Lemma \ref{s4_l2} is proved.
\end{proof}

\begin{remark}
In (\ref{s4_7}) there are no terms with boundary norms over the lateral
boundary (lines $x=0$ and $x=L$) because of the quasiperiodicity
condition.
\end{remark}
 From Lemma \ref{s4_l2} the basic result, Theorem \ref{s4_t1}, follows
immediately:
\begin{theorem}\label{s4_t1}
{\bf MRC-Modified Rayleigh Conjecture.} Fix $\ep>0$, however small,
and choose a positive integer $p$. Find
\begin{equation}\label{s4_9}
\min_{c_j}\|u_0+\sum_{|j|\leq p}c_j\p_j(x)v_j(y)\|:\,=m(p).
\end{equation}
Let $\{c_j(p)\}$ be the minimizer of (\ref{s4_9}). If $m(p)\leq\ep$, then
\begin{equation}\label{s4_10}
v(p)=\sum_{|j|\leq p}c_j(p)\p_j(x)v_j(y)
\end{equation}
satisfies the inequality
\begin{equation}\label{s4_11}
\||v-v(p)\||\leq c\ep,
\end{equation}
where $c=const>0$ does not depend on $\ep$. If $m(p)>\ep$, then there
exists $j=j(\ep)>p$ such that $m(j(\ep))<\ep$. Denote
$c_j(j(\ep)):\,=c_j(\ep)$ and $v(j(\ep)):\,=v_{\ep}$. Then
\begin{equation}\label{s4_12}
\||v-v_\ep\||\leq c\ep.
\end{equation}
\end{theorem}

\section{Numerical solution of the scattering problem}
According to the MRC method (Theorem \ref{s4_t1}), if the restriction of the 
incident field $-u_0(x,y)$ to $S_L$
 is approximated as in (\ref{s4_9}), then the series (\ref{s4_10})
approximates the scattered field in the entire region above the profile
$y=f(x)$. However, a numerical method that uses (\ref{s4_9}) does not
produce satisfactory results as
reported in \cite{petit} and elsewhere. Our own numerical
experiments confirm this observation. A way to overcome this difficulty
is to realize that the numerical approximation of the field 
$-u_0|_{S_L}$ can be carried out by using outgoing solutions 
described below.

Let $\bxi=(\xi_1,\xi_2)\in D_{-b}$, where $b>0$,
\[
D_{-b}:\,=\{(\xi_1,\xi_2)\ :\ -b\leq \xi_2\leq f(x), \quad 0\leq \xi_1\leq 
L\},
\]
and $g(\bx,\xi)$ be defined as in Section 3. Then $g(\bx,\xi)$ is an
outgoing solution satisfying $\Delta g+k^2g=0$ in $D_L$, according to (\ref{s3_3}).

To implement the MRC method numerically one proceeds as follows:
\begin{enumerate}
\item Choose the nodes $\bx_i,\ i=1,2,...,N$
 on the profile $S_L$. These points are used to approximate
$L^2$ norms on $S_L$.  
\item Choose points $\bxi^{(1)},\bxi^{(2)}, ..., \bxi^{(M)}$ in $D_{-b},\ 
M<N$.
\item Form the vectors ${\bf b}=(u_0(\bx_i))$, and 
${\bf a}^{(m)} = (g(\bx_i,\bxi^{(m)})),
\ i=1, 2,..., N, \ m= 1,2,...,M$. Let $\bf A$ be the $N\times M$ matrix containing
vectors ${\bf a}^{(m)}$ as its columns.
\item Find the Singular Value Decomposition of ${\bf A}$ . Use a predetermined $w_{min}>0$ to
eliminate its small singular values. Use the decomposition to compute
\[
r^{min}=\min\{\|{\bf b}+{\bf Ac}\|,\ {\bf c}\in \mathbb{C}^M\},
\]
where
\[
  \|{\bf a}\|^2=\frac 1{N}\sum_{i=1}^N |a_i|^2.
\]
\item {\bf Stopping criterion.} Let $\ep>0$. 
\begin{enumerate}
\item If $r^{min}\leq\ep$, then stop. Use the coefficients ${\bf c}=\{c_1,c_2,...,c_M\}$
 obtained
in the above minimization step to compute the scattered field by
\[
v(x,y)=\sum_{m=1}^M c_m g(x,y,\xi^{(m)}).
\]

\item If $r^{min}>\ep$, then increase $N, M$ by the order of 2, readjust
the location of points $\xi^{(m)}\in D_{-b}$ as needed, and repeat the procedure.
 \end{enumerate}
\end{enumerate}

We have conducted numerical experiments for four different profiles. In
each case we used $L=\pi, k=1.0$ and three values for the angle
$\theta$.  Table 1 shows the
resulting residuals $r^{min}$. Note that $\|{\bf b}\|=1$.
Thus, in all the considered cases, the MRC method
achieved $0.04\%$ to $2\%$ accuracy of the approximation.
Other parameters used in the experiments were chosen as follows:
$N=256,\ M=64,\ w_{min}=10^{-8},\ b=1.2$. The value of $b>0$, used in the
definition of $g$, was chosen experimentally, but the dependency of $r^{min}$ 
on $b$ was slight. The Singular Value Decomposition (SVD) is used in Step 4 since the vectors 
${\bf a}^{(m)},\ m=1,2,...,M$ may be nearly linearly dependent, which
leads to an instability in the determination of the minimizer $\bf c$.
According to the SVD method this instability is eliminated by cutting
off small singular values of the matrix $\bf A$, see e.g. \cite{numrec} for details. 
The cut-off value $w_{min}>0$ was chosen experimentally.
We used the truncated series (\ref{s3_01}) with $|j|\leq 120$
to compute functions $g(x,y,\xi)$. A typical run time on a 333 MHz PC was
about $40s$ for each experiment.

The following is a description of the profiles $y=f(x)$, the nodes
$\bx_i\in S_L$, and the poles $\xi^{(m)}\in D_{-b}$ used in  the computation of
$g(\bx_i,\xi^{(m)})$ in Step 3. For example, in profile I the $x$-coordinates of 
the $N$ nodes $\bx_i\in S_L$ are uniformly distributed on the interval
$0\leq x\leq L$. The poles $\xi^{(m)}\in D_{-b}$ were chosen as follows:
every fourth node $\bx_i$ was moved by a fixed amount $-0.1$ parallel to
the $y$ axis, so it would be within the region $D_{-b}$. The
location of the poles was chosen experimentally to give the smallest
value of the residual $r^{min}$.

{\bf Profile I.} $f(x)=sin(2x)$ for $0\leq x\leq L,\ t_i=iL/N,\ \bx_i=(t_i,f(t_i)),\ i=1,2,...,N,\ 
\xi^{(m)}=(x_{4m},y_{4m}-0.1),\ m=1,2,...,M$.

{\bf Profile II.} $f(x)=sin(0.2x)$ for $0\leq x\leq L,\ t_i=iL/N,\ \bx_i=(t_i,f(t_i)),\ i=1,2,...,N,\ 
\xi^{(m)}=(x_{4m},y_{4m}-0.1),\ m=1,2,...,M$.

{\bf Profile III.} $f(x)=x$ for $0\leq x\leq L/2$, $f(x)=L-x$ for
$L/2\leq x\leq L,\  t_i=iL/N,\ \bx_i=(t_i,f(t_i)),\ i=1,2,...,N,\ 
\xi^{(m)}=(x_{4m},y_{4m}-0.1),\ m=1,2,...,M$.

{\bf Profile IV.} $f(x)=x$ for 
$0\leq x\leq L,\ t_i=2iL/N,\ \bx_i=(t_i,f(t_i),\ i=1,...,N/2,
\  \bx_i=(L,f(2(i-N/2)L/N)),\ i=N/2+1,...,N,
\bxi^{(m)}=(x_{4m}-0.03,y_{4m}-0.05),\ m=1,2,...,M$.
In this profile $N/2$ nodes $\bx_i$ are uniformly distributed on its
slant part, and $N/2$ nodes are uniformly distributed on its vertical
portion $x=L$.

\begin{table}
\caption{Residuals attained in the numerical experiments.}

\begin{tabular}{c  c  c  }

\hline
Profile & $\theta$ & $r^{min}$ \\

\hline
I & $\pi/4$ & 0.000424  \\
  & $\pi/3$ & 0.000407  \\
  & $\pi/2$ & 0.000371  \\
\hline
II & $\pi/4$ & 0.001491  \\
   & $\pi/3$ & 0.001815  \\
   & $\pi/2$ & 0.002089  \\
\hline
III & $\pi/4$ & 0.009623  \\
   & $\pi/3$ & 0.011903  \\
   & $\pi/2$ & 0.013828  \\
\hline
IV & $\pi/4$ & 0.014398  \\
   & $\pi/3$ & 0.017648  \\
   & $\pi/2$ & 0.020451  \\
\hline
\end{tabular}

\end{table}

The experiments show that the MRC method provides a competitive
alternative to other methods for the computation of fields scattered
from periodic structures. It is fast and inexpensive. The results depend
on the number of the internal points $\xi^{(m)}$ and on their location. A
similar MRC method for the computation of fields scattered by a bounded
obstacle was presented in \cite{gutmanramm1}.


\begin{thebibliography}{99}

 \bibitem{alber} Alber, H.-D. A quasi-periodic boundary value problem for 
the Laplacian 
and the continuation of its resolvent.  Proc. Roy. Soc. Edinburgh Sect. A  
82  (1978/79), no. 3-4, 251--272.

 \bibitem{albertsen}
 Albertsen N.C., Chesneaux J.-M., Christiansen S., Wirgin A.,
  Comparison of four software packages applied to a scattering problem,
Mathematics and Computers in Simulation, 48, (1999), 307-317.


\bibitem{baran}  R. Barantsev,  
Concerning the Rayleigh hypothesis in the problem of scattering from 
finite bodies of arbitrary shapes,
Vestnik Lenungrad Univ., Math., Mech., Astron., 7 (1971) 56-62.

\bibitem{bon} A. Bonnet-Bendhia, Guided waves by eletromagnetic gratings 
and non-uniqueness examples for the diffraction problem,
Math. Math. in the Appl. Sci., 17, (1994), 305-338.

\bibitem{bonram}  A. Bonnet-Bendhia, K. Ramdani, Diffraction by an 
acoustic  grating perturbed by a bounded obstacle, Adv. Comp. Math., 16, 
(2002), 113-138.

 \bibitem{christiansen} S. Christiansen and R.E. Kleinman,
On a misconception involving point collocation and the Rayleigh
hypothesis,
IEEE Trans.Anten.Prop., 44,10, 1309-1316, 1996. 760

\bibitem{e} Eidus, D. M. Some boundary-value problems in infinite regions,
  Izv. Akad. Nauk SSSR Ser. Mat.  27, (1963) 1055--1080.

\bibitem{gutmanramm1}  Gutman, S., Ramm, A.G.,  
Numerical implementation of the MRC method for obstacle scattering problems,
J. Phys. A: Math. Gen., 35, (2002) 8065-8074.

\bibitem{k} Kato, T., Perturbation theory for linear operators,
Springer-Verlag, Berlin, 1995.

\bibitem{kazan1} Kazandjian L,
Rayleigh-Fourier and extinction theorem methods applied to scattering
and transmission at a rough solid-solid interface,
 J.Acoust.Soc.Am., 92,
1679-1691, 1992. 

\bibitem{kazan2} Kazandjian L,
Comments on "Reflection from a corrugated surface
revisited", [J. Acoust. Soc. Am., 96, 1116-1129 (1994)]" J. Acoust. Soc. Am.,
98, 1813-1814, (1995). 1245

\bibitem{millar1}  R. Millar, 
The Rayleigh hypothesis and a related least-squares solution 
to
the scattering problems for periodic surfaces and other scatterers,
Radio Sci., 8 (1973) 785-796.

\bibitem{millar2}  R. Millar, 
On the Rayleigh assumption in scattering by a periodic 
surface,
Proc. Camb. Phil. Soc., 69 (1971) 217-225.; 65 (1969) 773-791.

\bibitem{nazarov} Nazarov S., Plamenevskii B.,
 Elliptic problems 
in domains with piecewise smooth boundaries. de Gruyter Expositions in 
Mathematics, 13. Walter de Gruyter, Berlin, 1994. 

\bibitem{petit}  Petit R. (editor), Electromagnetic theory of gratings,
Topics in Current Physics, 22.
Springer-Verlag, Berlin-New York, 1980.

\bibitem{numrec} Press W.H., Teukolsky S.A., Vetterling W.T., Flannery B.P.
  [1992]
\textit{Numerical Recepies in FORTRAN,}
Second Ed., Cambridge University Press.

\bibitem{rammb1} A.G. Ramm, 
Scattering by Obstacles,
D. Reidel Publishing, Dordrecht, Holland, 1986.


\bibitem{r430} A.G. Ramm,
Modified Rayleigh Conjecture and Applications,
J. Phys. A: Math. Gen. 35 (2002) L357-L361.

\bibitem{r370} A.G. Ramm, G. Makrakis, Scattering by obstacles in acoustic 
waveguides, Spectral and scattering theory, 
in the book: editor A.G.RAMM, Plenum publishers, New York, 
1998, pp.89-110.

\bibitem{r189} A.G. Ramm,
Singularities of the inverses of Fredholm
 operators, Proc. of Roy. Soc. Edinburgh,
 102A, (1986), 117-121.

\bibitem{r3} A.G. Ramm,
Investigation of the scattering problem in
some
 domains with infinite boundaries I, II, Vestnik 7,
 (1963), 45-66; 19, (1963), 67-76.

\bibitem{r409} A.G. Ramm, M.Sammartino,
Existence and uniqueness of the scattering solutions
in the exterior of rough domains,

in the book "Operator Theory and Its Applications", Amer. Math.
Soc., Fields Institute

Communications vol.25, pp.457-472, Providence, RI, 2000.

(editors A.G.Ramm, P.N.Shivakumar, A.V.Strauss).

\bibitem{r475} A.G.Ramm, S.Gutman,
 {\it Modified Rayleigh Conjecture method for
multidimensional obstacle scattering problems}
(submitted)

\bibitem{r470}  A.G.Ramm, {\it Inverse Problems, vol. I,
II,}
Kluwer, Boston, 2004.

\bibitem{ral1} Rayleigh J.W., On the dynamical theory of gratings,
Proc. Roy. Soc. A, 79, (1907), 399-416.



\bibitem{voronovich} A. G. Voronovich, 
Wave Scattering from Rough Surfaces,
Springer, Berlin, 1996.




\end{thebibliography}
\end{document}